\documentclass [twoside, 11pt]{amsart}

\usepackage{amsmath}
\usepackage{amssymb}

\numberwithin{equation}{section} 

\newcommand{\cal}{\mathcal}

\newcommand{\amsc}{{\mathbb C}}  \newcommand{\amsz}{{\mathbb Z}}
\newcommand{\amsn}{{\mathbb N}}  \newcommand{\amsd}{{\mathbb D}}

\newcommand{\cbar}{\hat{{\mathbb C}} }

\newcommand{\lam}{\lambda}

\newcommand{\ph}{\varphi}
\newcommand{\al}{\alpha}

\newcommand{\jul}{{\cal J}_f}
\newcommand{\fat}{{\cal F}_f}

\newcommand{\conical}{\Lambda _c}

\newtheorem{theo}{Theorem}[section]

\newtheorem{coro}[theo]{Corollary}

\newtheorem{rema}[theo]{Remark}

\newtheorem{claim}[theo]{Claim}

\begin{document}

\title[]
{ \bf\Large {\huge T}he size of the Julia set of meromorphic
functions}

\date{\today}

\author[\sc Volker MAYER]{\sc Volker MAYER}
\address{Volker Mayer, Universit\'e de Lille I, UFR de Math\'ematiques, UMR 8524 du CNRS,
59655 Villeneuve d'Ascq Cedex, France}
\email{volker.mayer@univ-lille1.fr}

%

%

\begin{abstract}
We give a lower bound of the hyperbolic and the Hausdorff dimension
of the Julia set of meromorphic functions of finite order under very
general conditions.
\end{abstract}

\maketitle



\section{Introduction}
Let $f:\amsc\to \cbar$ be a meromorphic function. The Fatou set
$\fat$ is the set of points $z\in \amsc$ such that all the iterates
$f^n$ of $f$ are well defined and normal in some neighborhood of $z$.
The complement $\jul = \cbar \setminus \fat$ is the Julia set which
is the part of the sphere where the dynamics of the function are
complicated. However, some of the points $z\in \jul$ have the nice
property that a family of neighborhoods $\amsd(z,r_j)$, $r_j\to 0$,
can be zoomed conformally and with bounded distortion to a definite
size without going to infinity. These points built the so called set
of conical points $\conical $ (sometimes denoted by ${\cal J}_r$)
and its Hausdorff dimension is the hyperbolic dimension $HypDim(f)$
of the function $f$.

In this paper we show the following general theorem for a lower
estimate of the Hausdorff and hyperbolic dimension of the Julia set
of meromorphic functions of finite order. Concrete applications are
given in the last section.

\begin{theo}\label{theo main}
Let $f:\amsc \to \cbar$ be a meromorphic function of finite order
$\rho$ and suppose that:
\begin{itemize}
    \item[(i)] $f$ has (at least) one pole $b\in f(\amsc)$ which is
    not in the closure of the
    singular values $sing(f^{-1})$. Let $q $ be the
     multiplicity of the pole $b $.
    \item[(ii)] There is a neighborhood $D$ of $b$ and constants
    $K>0$, $\al_1 >-1-1/q$ such that
    $$ |f'(z) | \leq K |z| ^{\al_1} \quad for  \;\; z\in f^{-1}(D), \;
    |z| \to \infty \; .$$
\end{itemize}
Then
$$ Hdim (\jul ) \geq HypDim (f)  \geq \frac{\rho}{\al_1 +1 + 1/q}\;.$$
\end{theo}

It has been shown by G.M. Stallard \cite{st} that the Hausdorff and
the hyperbolic dimension of any meromorphic function is strictly
positif and also that there are transcendental meromorphic functions
with Julia set of Hausdorff dimension arbitrarily close to zero.

Our result is in the spirit of, and generalizes, work of J. Kotus
and M. Urba\'nski \cite{k,ku}. The main difference with \cite{k} is
that we need only a condition on one single pole of the function and
that we need not suppose that the post critical set stands away from
poles. Also, the condition (ii) is trivially satisfied for all
periodic functions (with $\al_1 =0$). Henceforth, our result applies
especially to all elliptic functions and in this case our theorem
reduces to the content of \cite{ku}. This and various other
applications are discussed in the last section of the paper.

Our result is also important for the new theory of thermodynamical
formalism for meromorphic functions of finite order in \cite{mu1,
mu2} since it yields an estimation for the topological pressure
which ensures the existence of geometric conformal measures (see in
particular Section 7 in \cite{mu1}).

We conclude this introduction in presenting now a first corollary of
our result.

\subsection{Exponential elliptic functions and continuity of Hausdorff dimension}

Compositions of exponential with elliptic functions are considered
in \cite{mu} and it is shown there that they all have hyperbolic
dimension two. These exponential elliptic functions do not fit
directly into our context. Indeed, our methods do not apply to
essential singularities. However they apply to
$$ f_d (z) = \lam\left( 1+\frac{f(z)}{d}\right)^d \quad ,\;z\in \amsc\; , $$
with $f$ any elliptic function, $d\geq 1$ entire and $\lam \in
\amsc\setminus \{0\}$. Such a function $f_d$ is of order $\rho =2$,
the number $\al_1$ is again zero and the maximal multiplicity of
poles is $dq$ with $q$ the maximal multiplicity of the poles of $f$.

\begin{coro}
The estimation $Hdim( {\cal J}_{f_d} )\geq HypDim (f_d)\geq
\frac{2dq}{dq+1}$ holds,
$$\lim_{d\to\infty}Hdim( {\cal J}_{f_d} )
= Hdim ( {\cal J}_{\lam \exp \circ f} ) =2 $$ and the same relation
is true for the hyperbolic dimension.
\end{coro}

\noindent

\section{Meromorphic Functions}
The reader may consult, for example, \cite{nev} or \cite{h} for a
detailed exposition on meromorphic functions and \cite{b} for their
dynamical aspects. We collect here the properties of interest for
our concerns.

Let $f:\amsc \to \cbar$ be a meromorphic function. We always suppose
that $f$ is of finite order $\rho=\rho(f)$. If $a\in \cbar$ is not
an omitted value, then we name the $a$--points $\{z_n(a)\; , \; n\in
\amsn\} = f^{-1} (a)$. The \it exponent of convergence \rm $ \rho_c
(f,a)$ of the series
$$ \Sigma (t,a)= \sum_n | z_n(a)|^{-t} $$
is defined by $\rho(f,a) = \inf \{t>0 \; ; \Sigma (t,a) < \infty
\}$. A Theorem of Borel (\cite{borel}, see also Corollary 2, p. 231
of \cite{h}) states that
$$ \rho(f,a) = \rho \quad \mbox{for all but at most two values} \;\; a\in
\cbar \; .$$

As usual we denote $\jul$ the Julia set and by $sing(f^{-1})$ the
set of singular values, i.e. critical and asymptotic values of $f$.
The following classes of functions are used (see \cite{b, el}):
\begin{itemize}
    \item $\cal S$ the functions $f$ with $sing(f^{-1})$ a
    finite set.
    \item $\cal B$ the functions with $sing(f^{-1})$ a
    bounded subset of the plane.
\end{itemize}
We also remark that a meromorphic function which is not entire has
either exactly one pole $b\in f^{-1}(\infty )$ which is an omitted
value (in which case the function is conjugate to a function of the
punctured plane) or some iterate of $f$ (actually $f^3$) has
infinitely many poles. The reason is Picard's Theorem.

\section{Iterated function systems and the Proof of Theorem \ref{theo main}}

Similar to \cite{ku, mu}, the main argument of the proof of Theorem
\ref{theo main} is that the function $f$ induces an iterated
function system for which we have a well developed theory thanks to
Mauldin - Urba\'nski \cite{mdu}.

Let $\Omega \subset \amsc $ be a simply connected domain with smooth
boundary and suppose that $\ph _j :\Omega \to \Omega _j=\ph_j(\Omega
) \subset \Omega $, $j\in \amsn $, are conformal mappings that are
uniformly contracting and which have uniformly bounded distortion
property. Suppose further that $\Omega_j \cap \Omega_i =\emptyset$
for all $i\neq j$. Then the limit set $L(S)$ of the iterated
function system $S=\{\ph_j:\Omega \to \Omega _j\, ; \, j\in \amsn
\}$ is defined by
$$ L(S) = \bigcap _{n=1}^\infty \bigcup_{j_1,...,j_n} \ph
_{j_1}\circ ...\circ \ph_{j_n} (\Omega ) \;\; .$$ Our proof of
Theorem \ref{theo main} relies on the following:\\

\noindent {\bf Theorem 3.15. of  \cite{mdu}:} For the Hausdorff
dimension of the limit set $L(S)$ we have the relation
$$Hdim ( L(S)) \geq \theta$$
where $\theta = \inf\{t>0\; ; \; \sum_{j\in \amsn} |\ph_j'(z_0)|^t
<\infty\}$, $z_0$ any point of $\Omega$.\\

\noindent {\bf Proof of theorem \ref{theo main}} Let $f$ be
meromorphic with pole $b$ of multiplicity $q$ and let $D$ be a
neighborhood of $b$ as in the conditions (i) and (ii) of Theorem
\ref{theo main}. Note that, near the pole $b$, the function $f$ has
the form $f(z) = \frac{g(z)}{(z-b)^q}$ with $g$ a holomorphic
function defined near $b$ and $g(b)\neq 0$.
 We may suppose that $D=\amsd(b,r)$ is a disk
centered at $b$ and with radius $r>0$ sufficiently small such that
firstly no singular value of $f$ belongs to $D^* =\amsd(b,2r)$ and
secondly
$$ |f'(z)| \simeq \frac{1}{|z-b|^{q+1}} \simeq |f(z)|^{1+1/q} \quad
for \;\; z\in D^*\setminus \{b\} \; .$$ Here and in the rest of the
paper $A\simeq B$ means that $\frac{A}{B}$ stands away from $0$ and
$\infty$ independently of the variables involved. We remark that the
last equation implies in particular that no critical point of $f$
belongs to $D^*$. Denote $V= f(D\setminus \{b\})$ and let $R>0$ such
that $\{|z|>R\}\subset V$.

Consider now the $b$--points $z_n=z_n(b) \in f^{-1}(b)\cap V$ and
denote $\ph_n$ the inverse branch of $f$ defined on $D$ by $\ph_n(b)
= z_n$.

\begin{claim}
There is $n_0\geq 1$ such that $\overline{\ph_n(D )} \subset U $ for
all $n \geq n_0$.
\end{claim}

\begin{proof}
Suppose $n\in\amsn$ such that there exists $y=\ph_n(x)\in
\overline{\ph_n(D )}\cap \{|z|=R\}$ and such that
$|z_n|=|\ph_n(b)|\geq 3R$. Denote $B_n=\ph_n(\amsd (x, r))$. By
Koebe's distortion theorem this set contains a disk $\amsd (y,
r|\ph_n'(x)|/4 )$. The same argument shows that
$$|\ph_n'(x)|\simeq |\ph_n'(b)| \simeq \frac{diam(\ph_n(D))}{diam
(D)} \geq \frac{|z_n|-R}{2r} \geq \frac{R}{r} \; .$$ Therefore,
there is $t=t(r,R)>0$ such that $B_n \supset \amsd (y, t)$ which
shows that $$\ph_n (D^*) \cap \{|z|=R\}$$ contains an arc of length
at least $t$. But this can happen only for finitely many
$n\in\amsn$, the sets $\ph_n (D^*)$, $n\in\amsn$, being disjoint.
The claim follows.
\end{proof}

Since $f:D\setminus \{b\} \to V$ is a unbranched covering there is,
for any $n\geq n_0$, an inverse branch $\psi_n$ of $f$ defined on
$\ph_n(D)$ with
$$ D_n = \Phi_n (D) = \psi_n \circ \ph_n (D) $$
compactly contained in $D$. Call $w_n = \psi _n (z_n)$. We now have
$${\cal S} = \left\{ \Phi_n : D \to D_n \; ; \; n\geq n_0 \right\}$$
a conformal iterated function system according to \cite{mdu} and it
remains to estimate the exponent of convergence of the associated
Poincar\'e series.

Suppose $w\in D$ such that $f^2(w)=b$ and set $z=f(w)\in V$. Then it
follows from condition (ii) that
$$|(f^2)'|(w)= |f'(w)||f'(z)|\lesssim |f(w)|^{1+1/q}|z|^{\al_1} =|z|^{\al_1 +1+1/q}.$$
Therefore
$$\Sigma_t =\sum _{n\geq n_0} |\Phi_n'(b)|^{t} = \sum _{n\geq n_0} |(f^{2})'(w_n)|^{-t} \gtrsim
\sum _{n\geq n_0} |z_n|^{-t(\al_1 +1+1/q)}\; .$$ Because of Koebe's
distortion theorem we may suppose that the pole $b$ is a Borel
point, meaning that the exponent of convergence of $\sum _n
|z_n|^{-t}$ is $\rho$ the order of $f$. We showed the following
estimation for the critical exponent $\theta$ of $\Sigma_t$:
$$\theta \geq \frac {\rho}{\al_1 +1+1/q}\; $$
and conclude that the limit set of the system $\cal S$ has Hausdorff
dimension at least $\theta$. The limit set of $\cal S$ being a
subset of the set of conical points of $f$ the Theorem is proven.
\hfill $\square$.

\begin{rema}
In many cases and especially in the examples given in the next
section the series $\Sigma_t$ is divergent at the critical exponent
$t=\theta$. The system $\cal S$ is then called hereditarily regular
and it is shown in \cite{mdu} that the estimation in Theorem
\ref{theo main} can be sharpened to
$$ Hdim (\jul ) \geq HypDim (f) > \frac{\rho}{\al _1 +1+1/q}.$$
\end{rema}

\section{Applications of Theorem \ref{theo main}}

We now discuss the conditions of Theorem \ref{theo main} and give
some applications and explicit examples. We will focus our attention
on the essential condition which is (ii). The first one (i) is often
automatically fulfilled. For example this is the case if $f\in {\cal
B}$, meaning that $sing(f^{-1})$ a bounded subset of the plane, and
if $f$ has infinity many poles. All the following examples are of
this kind.

\subsection{Periodic functions}

The polynomial growth condition (ii) on the derivative is always
satisfied for all periodic functions. Indeed, suppose $f:\amsc\to
\cbar$ is periodic, $b \in \amsc\setminus sing(f^{-1})$ is a pole
and $D=\amsd (b,r)$ a disk free of singular values. Then it follows
immediately from the periodicity that $|f'|$ is almost constant on
$f^{-1}(D)$. Therefore, condition (ii) is satisfied with $\al_1 =0$.

Consider, as a first example, the $1$--periodic family
 $$f_\lam (z) = \lam (tan \, z) ^m \;\; , \;\; m\in \amsn \;\; and \;\; \lam \in \amsc^* \; .$$
  All these
functions are of order $\rho =1$ and the (infinitely many) poles are
all of order $q=m$ which implies
$$Hdim({\cal J}_{f_\lam }) \geq HypDim (f_\lam )\geq
\frac{1}{1+1/m}= \frac{m}{m+1} \;\; for \; all \;\;\lam \in \amsc^*
\; .$$ This is precisely the lower bound found in \cite{k} and it is
shown there that, for this family, this bound is sharp.

Suppose now that $f:\amsc \to \cbar$ is
 an elliptic function, i.e. a meromorphic
function  for which there exists $\omega_1, \omega_2 \in \amsc$ with
$\Im \left( \frac{\omega_1}{\omega2}\right)
>0$ and such that
$$ f(a)=f(b) \;\; \mbox{if and only if} \;\; a=b+n\omega_1+m\omega_2$$
for some $n,m\in \amsz$. These functions are of order $\rho = 2$ and
 Theorem \ref{theo main} gives as the next Corollary the main result of \cite{ku}:

\begin{coro}
Suppose $f:\amsc\to \cbar$ is an elliptic function and $q$ the
maximal multiplicity of its poles. Then
$$ Hdim (\jul ) \geq HypDim (f)  \geq \frac{2}{1 + 1/q}=\frac{2q}{1 + q}$$
\end{coro}

 In the same way one can consider functions like $g(z) = f(z) +
f(\sqrt{2}z)$, $f$ elliptic, which are no longer periodic but the
same conclusion do hold. An other family of examples are functions
of the form $g=f\circ P$ with again $f:\amsc \to \cbar$ an elliptic
function and $P$ a polynomial of degree $d$. Then $g$ is of order
$\rho=2d$ and condition (ii) is satisfied with $\al_1 = d-1$. We
therefore get the estimation
$$ Hdim ({\cal J}_g) \geq HypDim (g) \geq \frac{2d}{d+1/p} \; .$$

\subsection{Solutions of differential equations}

Our result does also apply to solutions of suitable differential
equations. Let us illustrate this by inspecting the following
Ricatti equation
\begin{equation}\label{eq ricatti}
w'=P_0(z) +P_1(z) w + w^2
\end{equation}
were $P_0,P_1$ are polynomials of degree $d_0,d_1$ respectively.

Suppose $f:\amsc \to \cbar$ is a meromorphic solution of (\ref{eq
ricatti}) having infinitely many poles and at least one of them
outside $\overline{sing(f^{-1})}$. Then, by a result of Wittich (see
\cite[Thm 4.6.3]{h2}), $f$ is of finite order $\rho = 1
+\max\{\frac{1}{2}d_0,d_1\}$. Moreover, the growth condition (ii) of
Theorem \ref{theo main} is satisfied with $\al_1 = \max\{d_0,d_1\}$.
It easily follows then that
$$ Hdim(\jul) \geq HypDim(f) \geq \frac{\rho}{\al_1 +2} \geq
\frac{1}{2}$$ for any such function.

\subsection{Meromorphic functions with polynomial Schwarzian derivative}

The exponential and tangent functions are examples for which the
Schwarzian derivative $$ S(f) = \left( \frac{f''}{f'}\right)'
-\frac{1}{2} \left( \frac{f''}{f'}\right)^2$$ is constant. By
M\"obius invariance of $S(f)$, functions like
$$ \frac{e^z}{\lambda e^z +e^{-z}} \quad and \quad \frac{\lambda e^z}{e^z
-e^{-z}}$$ also have constant Schwarzian derivative. Examples for
which $S(f) $ is a polynomial are
$$ f(z) =\int_0^z \exp (Q(\xi )) \, d\xi  \quad , \quad Q\; a \;
polynomial,$$ and also
$$ f(z)= \frac{a\, Ai(z) +b\,Bi(z)}{c\, Ai(z) +d\,Bi(z)} \quad with
\quad ad-bc \neq0$$ and with $Ai$ and $Bi$ the Airy functions of the
first and second kind.

The asymptotic behavior of a general meromorphic solution $f$ of
$$ S(f) = P \quad , \quad P \;\; a \; polynomial,$$
 are well known and it turns out that such a
function satisfies the conditions of Theorem \ref{theo main} with
$\rho = d/2+1$ and $\al_1=d/2$ where $d=deg(P)$. Moreover, $f$ is of
divergence type and has only simple poles (details can be found in
Section 3 of \cite{mu2}). Consequently,
$$ Hdim (\jul ) \geq HypDim (f) > \frac{d+2}{d+4} \geq \frac{1}{2}.$$


\begin{thebibliography}{MM}


\bibitem[B]{b} W. Bergweiler, \em Iteration of meromorphic functions, \rm Bull. A.M.S. 29:2 (1993), 151-188.

\bibitem[Bo]{borel} \'E. Borel, \em Sur les z\'eros des fonctions enti\`eres,   \rm  Acta Math. 20 (1897), 357-396.

\bibitem[EL]{el} A.E. Eremenko and M.YU. Lyubich \em Dynamical properties of some classes of entire functions, \rm
Ann. Inst. Fourier, Grenoble 42, 4 (1992), 989-1020.

\bibitem[H1]{h} E. Hille, \em Analytic function theory, Vol II, \rm Ginn (1962).

\bibitem[H2]{h2} E. Hille, \em Ordinary differential equations in the complex domain, \rm Pure $\&$ Applied Mahtematics,
A Wiley--Interscience series of texts, Monographs $\&$ Tracts
(1976).

\bibitem[K]{k} J. Kotus, \em On the Hausdorff dimension of Julia sets
of meromorphic functions II, \rm Bull. Soc. Math. France 123 (1995),
33-46.

\bibitem[KU]{ku} J. Kotus and M. Urba\'nski, \em Hausdorff dimension and Hausdorff measures of Julia
sets of elliptic functions, \rm Bull. London Math. Soc. 35 (2003),
269-275.

\bibitem[MdU]{mdu} D. Mauldin and M. Urba\'nski, \em Dimensions and measures
in infinite iterated function systems, \rm Proc. London Math. Soc.
(3) 73 (1996) 105-154.



\bibitem[MyUr1]{mu} V. Mayer and M. Urba\'nski \em Exponential elliptic gives dimension two, \rm
Illinois J. Math., \rm to appear.


\bibitem[MyUr2]{mu1} V. Mayer and M. Urba\'nski, \em Geometric
thermodynamical formalism and real analyticity for meromorphic
functions of finite order, \rm preprint.

\bibitem[MyUr3]{mu2} V. Mayer and M. Urba\'nski, \em Thermodynamical
formalism and multifractal analysis  for meromorphic functions of
finite order, \rm preprint.



\bibitem[Nev]{nev} R. Nevanlinna, \em Analytic functions, \rm Springer Verlag, Berlin (1970).


\bibitem[St]{st} G.M. Stallard, \em The Hausdorff dimension of Julia sets of meromorphic functions,
\rm J. London Math. Soc. (2) 49 (1994), 281- 295.


\end{thebibliography}
 \end{document}